\documentclass{article}

\usepackage{nips15submit_e,times}
\usepackage{npbayes}
\usepackage{verbatim}
\usepackage{graphicx}
\usepackage{amsmath, amssymb, amsthm}
\usepackage{url}

\title{Edge-exchangeable graphs and sparsity}

\author{Tamara Broderick \\ Department of EECS \\ Massachusetts Institute of Technology \\
    \texttt{tbroderick@csail.mit.edu}
\And Diana Cai \\ Department of Statistics \\ University of Chicago \\
    \texttt{dcai@uchicago.edu}}

\begin{document}

\nipsfinalcopy
\maketitle

\begin{abstract}
A known failing of many popular random graph models is that the Aldous--Hoover
Theorem guarantees these graphs are dense with probability one; that is, the
number of edges grows quadratically with the number of nodes. This behavior is
considered unrealistic in observed graphs. We define a notion of edge
exchangeability for random graphs in contrast to the established notion of
infinite exchangeability for random graphs---which has traditionally relied on
exchangeability of nodes (rather than edges) in a graph. We show that, unlike
node exchangeability, edge exchangeability encompasses models that are known to
provide a projective sequence of random graphs that circumvent the
Aldous--Hoover Theorem and exhibit sparsity, i.e., sub-quadratic growth of the
number of edges with the number of nodes. We show how edge-exchangeability of
graphs relates naturally to existing notions of exchangeability from clustering
(a.k.a.\ partitions) and other familiar combinatorial structures.
\end{abstract}

\section{Introduction}

As graph, or network, data become more ubiquitous and large-scale---in the form
of social networks, collaboration networks, networks expressing biological
interactions, etc.---a number of probabilistic models for graphs have been proposed.
However, it has become apparent that many of the most popular probabilistic models
for graphs are fundamentally misspecified in a way that worsens as the graphs scale
to larger numbers of nodes \citep{orbanz:2015:bayesian}. At the heart of the
problem is that generative probabilistic modeling relies on what seems at first
to be a very weak assumption, that of \emph{exchangeability}. Exchangeability
is essentially the idea that seeing our data in a different order does not change
its distribution and is a much weaker assumption than the popular ``independent
and identical distribution (iid)'' assumption. In graphs, this assumption has
historically taken the following form: we assume that if we relabeled our nodes,
it would not change the probability of the graph
\citep{aldous:1985:exchangeability,Hoover79,orbanz:2015:bayesian}.
Exchangeability assumptions are fundamentally tied with probabilistic modeling since they
imply the existence of parameters, likelihoods, and priors by \emph{de Finetti's Theorem}.
The particular version of de Finetti's Theorem for graphs with this form of
\emph{node-exchangeability} is known as the \emph{Aldous--Hoover Theorem}
\cite{aldous:1985:exchangeability,Hoover79}.
Notably, the Aldous-Hoover Theorem for exchangeable graphs implies that,
if we assume our graph is node-exchangeable, it must be that our graph is
\emph{dense} \cite{orbanz:2015:bayesian}.
That is, the number of edges in a dense graph grows quadratically with the number of nodes.

But most real-world graphs have been observed to be \emph{sparse}
\citep{newman2005power};
that is, the number of edges grows sub-quadratically in the number of nodes.
While there are ad hoc solutions to this mismatch between model and data, they
have other undesirable properties, such as lacking \emph{projectivity}---a
property that facilitates handling streaming data, performing distributed data
analysis, and consistent hierarchical modeling. \citet{caron:2015:bayesian} have
recently suggested an example model that has some desirable sparse scaling
properties. They consider an alternative form of exchangeability in the sense of
independent increments of subordinators. By contrast, we here consider a new
form of exchangeability for graphs where we consider permuting the \emph{edges}
rather than the nodes. In the remainder of the current work, we both describe
node-exchangeability in more detail and introduce our new concept of
\emph{edge-exchangeability}. We show how edge-exchangeability of graphs relates
naturally to existing notions of exchangeability from clustering
(a.k.a.\ partitions), feature allocations \cite{broderick:2013:cluster,broderick:2013:feature}, and other
combinatorial structures \cite{broderick:2015:combinatorial,broderick:2014:posteriors}.
We describe how the \citet{caron:2015:bayesian} model fits into our framework.
We outline remaining connections and characterizations to be made in future work.

\section{Edge-exchangeability in graphs}
\label{sec:edge-exchangeable}

An undirected graph is defined by a set of nodes (i.e., vertices) called $V$ and a set of edges $E$. In particular, each element $e$ of $E$ is an unordered set $\{u,v\}$ of two nodes $u, v \in V$ and represents a link between $u$ and $v$.
We consider $V$ to be a set where each element is unique, but we allow $E$ to be a multiset; that is, we allow edges to potentially occur multiple times in $E$.
We define \emph{active nodes} to be those nodes that appear in some edge, where an edge between
a vertex and itself counts for this purpose.
In what follows, we take the approach that we are only interested in active nodes. In this case,
an undirected graph may be characterized
by specifying only its edge set $E$. That is,
we can consider a graph as an unordered set of tuples, which represent the edges of the
graph. In this case, we can obtain the (active) node set by
representing all the unique elements of $E$: $V = \bigcup_{e \in E} e$.

\begin{example} \label{ex:active}
Consider the graph (containing only active nodes) defined by $E = \{\{2,3\},\{1,4\},\{3,6\}, \{6,6\},\{3,6\}\}$. 
Then the active nodes are $V = \{1,2,3,4,6\}$.
It is not important that ``$5$'' does not appear in $V$ since the elements of $V$ are treated as arbitrary labels.
\end{example}

To consider exchangeability in graphs, it has been traditional to think
of a sequence of (random) graphs $E_1, E_2, \ldots$ such that $E_m \subseteq E_n$ for any
$m < n$ and such that $E_n$ represents the edges between the vertices with
labels in $[n] := \{1,\ldots,n\}$. One typically thinks of the step of going from
$E_{n}$ to $E_{n+1}$ as adding in a new vertex labeled $n+1$ (and all its edges
to existing vertices and itself). Then we say that the graphs are
\emph{infinitely exchangeable} if, for any positive integer $n$, we have that
permuting the vertex labels in $[n]$ does not change the probability of these
graphs. This form of exchangeability, which we will refer to as
\emph{node exchangeability} implies, via the Aldous--Hoover Theorem \cite{Hoover79, aldous:1985:exchangeability}, that our graphs must be dense \cite{orbanz:2015:bayesian}. Thus, if we wish to maintain some notion of exchangeability while modeling sparse graphs, it behooves us to consider alternative forms of exchangeability.

\begin{example} \label{ex:traditional}
An example realization of a sequence of random graphs when considering infinite exchangeability
in the traditional sense is given by
$$
	E_1 = \emptyset,
	E_2 = \{\{1,2\}\},
	E_3 = \{\{1,2\}\},
	E_4 = \{\{1,2\}, \{1,4\},\{2,4\},\{3,4\}\},
	\ldots.
$$
Note that it is not necessary to specify the node set $V_1$, $V_2$, $V_3$, $\ldots$ here since
it is understood that $V_{n} = [n]$. Further note that $V_n$ here represents all the nodes, not only the
active ones.
\end{example}

In the present work, we instead introduce a notion of exchangeability of the \emph{edges} rather than the \emph{vertices}. In particular, consider a new sequence of graphs $E_1, E_2, \ldots$, where $E_m \subseteq E_n$ for any $m < n$. Now we think of $E_{n+1}$ as adding some new edges relative to $E_{n}$, but these new edges need not be connected to any particular vertex. We can make the step on which we add in edge $e$ explicit by augmenting the edge set $E$. In particular, define a \emph{step-augmented graph} $E'$ as a collection of tuples, where the first element is the edge and the second element is the step on which the edge is added. 

\begin{example}
There is exactly one sequence of step-augmented graphs that corresponds to the graph sequence in \ex{traditional}. It is
\begin{align*}
	E'_1 &= \emptyset, \\
	E'_2 &= \{(\{1,2\},2)\}, \\
	E'_3 &= \{(\{1,2\},2)\}, \\
	E'_4 &= \{(\{1,2\},2), (\{1,4\},4),(\{2,4\},4),(\{3,4\},4)\}.
\end{align*}
In the traditional exchangeability setup of \ex{traditional}, the step of an edge is always the maximum node value in that edge.
\end{example}

\begin{example} \label{ex:augm}
In our new setup, the step need not be the maximum node value. For instance, consider a step-augmentation of \ex{active}:
$$
	E'_4 = \{(\{2,3\},1),(\{1,4\},4),(\{3,6\},1), (\{6,6\},3), (\{3,6\},3)\}.
$$
This augmentation would be equivalent to the graph sequence
\begin{align*}
	E_1 &= \{\{2,3\}, \{3,6\}\}, \\
	E_2 &= \{\{2,3\}, \{3,6\}\}, \\
	E_3 &= \{\{2,3\}, \{3,6\}, \{6,6\}\}, \{3,6\}\}, \\
	E_4 &= \{\{2,3\},\{1,4\},\{3,6\}, \{6,6\}, \{3,6\} \}.
\end{align*}
\end{example}

Now suppose we treat the \emph{edges} (rather than the nodes) in this sequence of edge sets as arbitrary labels; if two edges are the same, they have the same labels, and otherwise they have different labels. This representation is technically a superclass of the objects defined above since there is additional structure in the edge labels of graphs; graph edge labels may agree in one element but not the other.
\begin{example} \label{ex:labels}
Using an order of appearance scheme \cite{broderick:2013:feature} to index the labels, $E'_4$ in \ex{augm} becomes $\{(\phi_1,1), (\phi_2,1), (\phi_3,3), (\phi_1,3), (\phi_4,4)\}$.
\end{example}
We note that if we consider only the uniqueness of the labels in $E'_n$ and not their actual values, the information in this structure can be expressed as a set of subsets of the steps $[n]$. That is, for each unique label $\phi$ that occurs in any tuple in $E'_n$, we let $A_{\phi}$ be the collection of steps that co-occur in a tuple with $\phi$. In \ex{labels}, $A_{\phi_1} = \{1,3\}$. Then let $C_n$ be the collection of $A_{\phi}$ for all values of $\phi$. Thus, $C_n$ is a set of subsets of $[n]$ since $E'_n$ contains only edges added on step $n$ or earlier. We call the sequence $(C_n)$ the \emph{step collection sequence} of a sequence of graphs.
\begin{example}
The $C_4$ corresponding to $E'_4$ in \ex{labels} is $\{\{1,3\}, \{1\}, \{3\}, \{4\}\}$.
\end{example}
Finally, now, we recognize the step collection $C_n$ as a familiar combinatorial object: either a partition (a.k.a.\ clustering) \cite{pitman:1995:exchangeable}, feature allocation \cite{broderick:2013:feature}, or a trait allocation (defined below and reminiscent of \cite{broderick:2014:posteriors}); we explain each of these connections below and show how they give a natural notion of exchangeability in the \emph{edges} of a graph.

\subsection{Partition connection} \label{sec:part}

First consider the connection to partitions.
In this case, suppose that each index in $[n]$ appears exactly once across all
of the subsets of $C_n$. This assumption on $C_n$ is equivalent to assuming that in the original
graph sequence $E_1, E_2, \ldots$, we have that $E_{n+1}$ always has exactly one
more edge than $E_{n}$. In this case, $C_n$ is exactly a
\emph{partition} of $[n]$; that is, $C_n$ is a set of mutually exclusive and
exhaustive subsets of $[n]$. If the edge sequence $(E_n)$ is random,
then $(C_n)$ is random as well.

We say that a partition sequence $C_1, C_2, \ldots$,
where $C_n$ is a (random) partition of $[n]$ and $C_m \subseteq C_n$ 
for all $m \le n$, is infinitely exchangeable if, for all $n$,
permuting the indices of $n$ does not change the distribution
of the (random) partitions \cite{pitman:1995:exchangeable}. Permuting the indices $[n]$
in the partition sequence $(C_m)$
corresponds to permuting the order in which edges are added in our graph
sequence $(E_m)$. Contrast this with the traditional form of
exchangeability in graphs (node exchangeability) as described above.

Recall further that the \emph{Kingman paintbox theorem} \citep{kingman:1978:representation} tells us that we have an infinitely exchangeable partition sequence if and only if we can find a sequence of (potentially random) probabilities $p_1, p_2, \ldots$ such that $p_k \in (0,1)$ and $\sum_{k=1}^{\infty} p_k = 1$ and such that drawing partition elements according to these probabilities yields the same partition distribution as our original random partition. The sequence $(p_k)_{k=1}^{\infty}$ is called the \emph{Kingman paintbox}. In the graph domain, we can interpret the Kingman paintbox probability $p_k$ as the probability of a particular edge in the graph.

\begin{example} \label{ex:graph_part}
We consider a generative model proposed by \citet{caron:2015:bayesian}. Let $W =
\sum_{k=1}^{\infty} w_k \delta_{\phi_k}$ be a random measure such that the pairs
$\{(w_k,\phi_k)\}_{k=1}^{\infty}$ are generated from a Poisson point process
with rate measure $\nu(dw, d\phi) = \nu(dw) G(d\phi)$ for some proper
distribution $G$. We assume that $\nu$ is a positive measure with support on $\mathbb{R}_{+}$. In this case, $W$
is a \emph{completely random measure}. For $n=1,2,\ldots$, we draw whether the
graph acquires edge $\{i,j\}$ at step $n$ according to the distribution
$(p_{\{i,j\}})_{i,j}$ where
$$
	p_{\{i,j\}} \propto \left\{ \begin{array}{ll}
		2 w_{i} w_{j} & i\ne j \\
		w_{i}^{2} & i = j
		\end{array} \right. .
$$
Defining this distribution requires that the sum over the $w_{i}
w_{j}$ be finite, but this finiteness condition will hold for a wide range of models, including the
generalized gamma process \cite{hougaard1986survival, klein2003survival} with $\phi \in [0,\alpha]$ for some fixed $\alpha > 0$, as considered by \citep{caron:2015:bayesian}. In fact, this construction returns exactly the generative model of \citet{caron:2015:bayesian} if the total number of edges $N$ is chosen to be Poisson with appropriate rate parameter.\footnote{This observation makes use of the fact that a gamma process prior paired with a Poisson likelihood process yields the same distribution as a Dirichlet process paired with a multinomial likelihood and a Poisson-distributed number of data points.}
\end{example}

\subsection{Feature allocation connection} \label{sec:feat}

Next we notice that it need not be the case that exactly one edge is added at
each step of the graph sequence, e.g.\ between $E_{n}$ and $E_{n+1}$. If we
allow multiple unique edges at any step, then the step collection $C_n$ is just a set
of subsets of $[n]$, where each subset has at most one of each index in $[n]$.
Suppose that any $m$ belongs to only finitely many subsets in $C_n$ for any $n$.
That is, we suppose that only finitely many edges are added to the graph at any
step. Then $C_n$ is an example of a \emph{feature allocation}
\citep{broderick:2013:feature}. Again, if $(E_n)$ is random, then $(C_n)$ is random as well.

We say that a (random) feature allocation sequence $(C_m)$ is 
infinitely exchangeable if, for any $n$, permuting the indices of $[n]$ does not
change the distribution of the (random) feature allocations
\cite{broderick:2013:cluster, broderick:2013:feature}.
In this feature allocation case, permuting the indices $[n]$ in the sequence $(C_m)$ corresponds to
permuting the steps when edges are added in the edge sequence $(E_m)$.

Similarly to the partition case in \mysec{part}, we can apply known results from feature allocations to characterize edge-exchangeable graph models of this form. For instance, we know that the \emph{feature paintbox} \cite{broderick:2013:feature} characterizes distributions over feature allocations (and therefore this sequence of edge-exchangeable graphs) just as the Kingman paintbox characterizes distributions over partitions (and therefore the edge-exchangeable graphs with exactly one new edge per step).

We may also consider feature paintbox distributions with extra structure. For instance, we may say that an edge-exchangeable graph sequence (with multiple unique edges per step) has an \emph{edge-exchangeable graph probability function} (EGPF) if the probability of the graph can be expressed as a function only of the total number of steps $N$ and the edge multiplicities (and where the probability is symmetric in the edge multiplicities). This definition directly corresponds to the notion of an exchangeable feature probability function (EFPF) \cite{broderick:2013:feature,broderick:2013:cluster} on the feature allocations $(C_n)$.\footnote{Both functions are related to exchangeable partition probability functions (EPPFs) \citep{pitman:1995:exchangeable} for partitions.}

Also, we may define a \emph{graph frequency model} as built around a random measure $B = \sum_{k=1}^{\infty} V_k \delta_{\phi_k}$. In particular, we draw a random graph conditional on $B$ as follows. For each step index $n$, independently make a Bernoulli draw with success probability $V_k$. If the draw is a success, the edge indexed by $k$ appears at time $n$. Otherwise it does not.

Then Theorem 17 (``Equivalence of EFPFs and feature frequency models'') from
\cite{broderick:2013:feature} translates into a theorem about edge-exchangeable graphs as follows.
\begin{theorem} \label{thm:egpf}
Let $\lambda$ be a non-negative random variable (which may have some arbitrary joint law with the frequencies in a graph frequency model). We can obtain an edge-exchangeable graph by generating an edge-exchangeable graph from a graph frequency model and then, for each time $n$, including an independent $\pois(\lambda)$-distributed number of unique edges (which are different from those previously generated and which will never appear again). A graph of this type has an EGPF. Conversely, every graph with an EGPF has the same distribution as one generated by this construction for some joint distribution of $\lambda$ and the edge frequencies.
\end{theorem}

\begin{example}
We consider a graph frequency model related to the model from \ex{graph_part}. In particular, let $W = \sum_{k=1}^{\infty} w_k \delta_{\phi_k}$ be a completely random measure with rate measure $\nu(dw, d\phi) = \nu(dw) G(d\phi)$ for some proper distribution $G$. We now assume that $\nu$ is a positive measure with support on $[0,1]$. For $n=1,2,\ldots$, we draw whether the graph has an edge $\{i,j\}$ at time step $n$ as
$$
	\bern(q_{\{i,j\}}) \quad
		\textrm{ where } \quad
		q_{\{i,j\}} := \left\{ \begin{array}{ll}
			2 w_i w_j & i \ne j \\
			w_i^2 & i = j
			\end{array} \right. .
$$
Since the graph frequency representation is given explicitly, the EGPF existence follows by \thm{egpf}.
\end{example}

\subsection{Trait allocation connection}

Finally, we may consider the case where at every step, any non-negative (finite)
number of edges may be added \emph{and} those edges may have non-trivial
(finite) multiplicity; that is, the multiplicity of any edge at any step can be any
non-negative integer. By contrast, in \mysec{feat}, each unique edge occurred at most
once at each step. In this case, the step collection $C_n$ is a set of subsets of $[n]$.
The subsets need not be unique or exclusive since we assume any number of edges may be
added at any step. And the subsets themselves are multi-sets since an edge
may be added with some multiplicity at step $n$. We say that $C_n$ is a 
\emph{trait allocation},\footnote{We
introduce this terminology following the discussion in \cite{broderick:2014:posteriors}.}
which we define as a generalization of a feature allocation
where the subsets of $C_n$ are multisets. As above, if $(E_n)$ is random, $(C_n)$ is as well.

We say that a (random) trait allocation sequence $(C_m)$ is infinitely exchangeable
if, for any $n$, permuting the indices of $[n]$ does not change the distribution
of the (random) trait allocation. Here, permuting the indices of $[n]$ 
corresponds to permuting the steps when edges are added in the edge
sequence $(E_m)$.

\begin{example}
Consider again the construction from \ex{graph_part}. Now suppose that at each step, we add $n_{\{i,j\}}$ instances of edge $\{i,j\}$, where $n_{\{i,j\}}$ is drawn from some distribution $h(\cdot | \theta_{\{i,j\}})$, where $\theta_{\{i,j\}}$ equals $2 w_i w_j$ for $i \ne j$ and $w_i^2$ for $i = j$. The case where $h(\cdot | \theta)$ is Poisson with mean parameter $\theta$ and we take exactly one step can be seen as another interpretation of the model in \citet{caron:2015:bayesian}.
\end{example}

\section{Conclusion}
\label{sec:concl}

In this work, we have defined a notion of edge exchangeability for random graphs in contrast to the 
traditional notion of node exchangeability for random graphs. 
While the Aldous-Hoover Theorem guarantees that node-exchangeable random graphs
must be dense with probability one, 
we have seen in \ex{graph_part} that edge-exchangeability
encompasses models that are known to provide a projective sequence of
random graphs that circumvent this theorem and exhibit sparsity.
It remains to more fully characterize the asymptotic properties of edge-exchangeable
random graphs. For one, we have considered how certain types of structure (e.g., the EGPF
in \mysec{feat}) affect edge-exchangeable random graphs (\thm{egpf}). 
But perhaps an even more natural type of structure for graphs would be the 
idea that the probability of a graph depends only on the number of times a node
occurs. We believe that this will yield a probability function structure like the EGPF.
Moreover, it remains to consider power laws and other asymptotic behaviors in 
graphs in the style of \cite{gnedin:2007:notes} for partitions and
\cite{broderick:2012:beta} for feature allocations.

\bibliographystyle{plainnat}
\bibliography{bnp}

\begin{thebibliography}{15}
\providecommand{\natexlab}[1]{#1}
\providecommand{\url}[1]{\texttt{#1}}
\expandafter\ifx\csname urlstyle\endcsname\relax
  \providecommand{\doi}[1]{doi: #1}\else
  \providecommand{\doi}{doi: \begingroup \urlstyle{rm}\Url}\fi

\bibitem[Aldous(1985)]{aldous:1985:exchangeability}
D.~J. Aldous.
\newblock \emph{Exchangeability and related topics}.
\newblock Springer Lecture Notes in Math, 1985.

\bibitem[Broderick et~al.(2012)Broderick, Jordan, and
  Pitman]{broderick:2012:beta}
T.~Broderick, M.~I. Jordan, and J.~Pitman.
\newblock Beta processes, stick-breaking, and power laws.
\newblock \emph{Bayesian Analysis}, 7\penalty0 (2):\penalty0 439--476, 2012.

\bibitem[Broderick et~al.(2013{\natexlab{a}})Broderick, Jordan, and
  Pitman]{broderick:2013:cluster}
T.~Broderick, M.~I. Jordan, and J.~Pitman.
\newblock Cluster and feature modeling from combinatorial stochastic processes.
\newblock \emph{Statistical Science}, 2013{\natexlab{a}}.

\bibitem[Broderick et~al.(2013{\natexlab{b}})Broderick, Pitman, and
  Jordan]{broderick:2013:feature}
T.~Broderick, J.~Pitman, and M.~I. Jordan.
\newblock Feature allocations, probability functions, and paintboxes.
\newblock \emph{Bayesian Analysis}, 8\penalty0 (4):\penalty0 801--836,
  2013{\natexlab{b}}.

\bibitem[Broderick et~al.(2014)Broderick, Wilson, and
  Jordan]{broderick:2014:posteriors}
T.~Broderick, A.~C. Wilson, and M.~I. Jordan.
\newblock Posteriors, conjugacy, and exponential families for completely random
  measures.
\newblock \emph{arXiv preprint arXiv:1410.6843}, 2014.

\bibitem[Broderick et~al.(2015)Broderick, Mackey, Paisley, and
  Jordan]{broderick:2015:combinatorial}
T.~Broderick, L.~Mackey, J.~Paisley, and M.~I. Jordan.
\newblock Combinatorial clustering and the beta negative binomial process.
\newblock \emph{IEEE TPAMI}, 2015.

\bibitem[Caron and Fox(2015)]{caron:2015:bayesian}
F.~Caron and E.~B. Fox.
\newblock Bayesian nonparametric models of sparse and exchangeable random
  graphs.
\newblock \emph{arXiv preprint arXiv:1401.1137v3}, 2015.

\bibitem[Gnedin et~al.(2007)Gnedin, Hansen, and Pitman]{gnedin:2007:notes}
A.~Gnedin, B.~Hansen, and J.~Pitman.
\newblock Notes on the occupancy problem with infinitely many boxes: general
  asymptotics and power laws.
\newblock \emph{Probab. Surv.}, 4:\penalty0 146--171, 2007.

\bibitem[Hoover(1979)]{Hoover79}
D.~N. Hoover.
\newblock Relations on probability spaces and arrays of random variables.
\newblock Preprint, Institute for Advanced Study, Princeton, NJ, 1979.

\bibitem[Hougaard(1986)]{hougaard1986survival}
P.~Hougaard.
\newblock Survival models for heterogeneous populations derived from stable
  distributions.
\newblock \emph{Biometrika}, 73\penalty0 (2):\penalty0 387--396, 1986.

\bibitem[Kingman(1978)]{kingman:1978:representation}
J.~F.~C. Kingman.
\newblock The representation of partition structures.
\newblock \emph{Journal of the London Mathematical Society}, 2\penalty0
  (2):\penalty0 374--380, 1978.

\bibitem[Klein and Moeschberger(2003)]{klein2003survival}
J.~P. Klein and M.~L. Moeschberger.
\newblock \emph{Survival analysis: techniques for censored and truncated data}.
\newblock Springer Science \& Business Media, 2003.

\bibitem[Newman(2005)]{newman2005power}
M.~E.~J. Newman.
\newblock Power laws, {P}areto distributions and {Z}ipf's law.
\newblock \emph{Contemporary physics}, 46\penalty0 (5):\penalty0 323--351,
  2005.

\bibitem[Orbanz and Roy(2015)]{orbanz:2015:bayesian}
P.~Orbanz and D.~M. Roy.
\newblock Bayesian models of graphs, arrays and other exchangeable random
  structures.
\newblock \emph{Pattern Analysis and Machine Intelligence, IEEE Transactions
  on}, 37\penalty0 (2):\penalty0 437--461, 2015.

\bibitem[Pitman(1995)]{pitman:1995:exchangeable}
J.~Pitman.
\newblock Exchangeable and partially exchangeable random partitions.
\newblock \emph{Probability theory and related fields}, 102\penalty0
  (2):\penalty0 145--158, 1995.

\end{thebibliography}

\end{document}